\documentclass{glasnik}
\usepackage{glasnik}
\usepackage{graphicx}
\usepackage{url}

\newcommand{\Z}{\mathbb{Z}} 
\newcommand{\N}{\mathbb{N}} 

\begin{document}

\title[On $k$-th power Diophantine triples]{On $k$-th power Diophantine triples of the form $\{a^k, b, c\}$}
\UniCountry{University of Salzburg, Austria}
\author[C. Fuchs and M. Sch\"{o}nauer]{Clemens Fuchs and Miriam Sch\"{o}nauer}
\address[C. Fuchs]{Mathematics Department\\University of Salzburg\\ 5020 Salzburg\\ Austria}
\email{clemens.fuchs@plus.ac.at}
\address[M. Sch\"{o}nauer]{Mathematics Department\\University of Salzburg\\ 5020 Salzburg\\ Austria}
\email{miriam.schoenauer@plus.ac.at}



\abstract{In this paper, we prove that there are no $k$-th power Diophantine triples of the form $\{a^k,b,c\}$ for $k\geq 3$ and $1<a^k<b<c$.}

\maketitle

\section{Introduction}
\label{intro}

  The mathematician Diophantus of Alexandria was the first to consider the problem of finding a set of four rational numbers such that the product of any two elements increased by one is a perfect square. 
  He discovered a set of numbers that satisfy this property: $\left\{\frac{1}{16},\frac{33}{16},\frac{17}{4},\frac{105}{16}\right\}$.
  The first set of integers with this property was found by Fermat, namely $\left\{1,3,8,120\right\}$.
  Now, we call a set of $m$ positive integers $\{a_1,\dots,a_m\}$ a Diophantine $m$-tuple if $a_ia_j+1$ is a perfect square for all $1\leq i < j\leq m$.
  
  The question of how large these sets can be has created a wide and active field of research. To name only a few, one of the first important results was by Baker and Davenport \cite{BakerDavenport69}, who showed that for a Diophantine quadruple of the form $\{1,3,8,d\}$ the positive integer $d$ must be $120$. Thus, showing that this quadruple could not be extended to a quintuple.
  
  The first absolute bound for Diophantine $m$-tuples was established in 2001 by Dujella \cite{Dujella01}, who showed that $m\leq 8$. In 2004 the same author further sharpened the upper bound and proved that there are no Diophantine sextuples and only finitely many quintuples, see \cite{Dujella04}. Finally, He, Togbé and Ziegler \cite{HeTogbeZiegler19} proved that there are no Diophantine quintuples. A detailed overview of the literature and results on Diophantine $m$-tuple can be found on the webpage \cite{DujellaWeb}, see also \cite{Dujella24}.

  As a generalization of Diophantine $m$-tuples one may consider sets of $m$ positive integers $\{a_1,\dots, a_m\}$ such that $a_ia_j+1$ is a perfect $k$-th power for all $1\leq i < j\leq m$ and $k\geq 3$. Such sets are called $k$-th power Diophantine $m$-tuples. Examples of such sets for $k = 3$ and $k = 4$ are $\{2, 171, 25326\}$ and $\{1352, 9539880, 9768370\}$, respectively. Bugeaud and Dujella \cite{BugeaudDujella03} examined $k$-th power Diophantine $m$-tuples in 2003 and established absolute bounds for them, namely if
  $$
    C_k:=\sup\{|\mathit{A}|: \mathit{A} \text{ is a } k \text{-th power Diophantine } m \text{-tuple}\}
  $$
  for $k\geq 3$ then it holds that $C_3\leq 7, C_4\leq 5$, $C_5\leq 5$, $C_k \leq 4$ for $6\leq k\leq 176$ and $C_k\leq 3$ for $k\geq 177$. As pointed out by Bérczes, Dujella, Hajdu and Luca \cite{Berczes11} a slight inaccuracy in the proof of \cite[Corollary 4]{BugeaudDujella03} resulted in the bound $C_5\leq 4$ instead of $C_5\leq 5$.  
  Additionally, it follows from \cite{HeTogbeZiegler19} that for even $k$, it holds that $C_k\leq 4$.
  In 2004 Bugeaud \cite{Bugeaud04} showed that for a $k$-th power Diophantine triple of the form $\{1,a,b\}$, $k$ may not exceed $74$. In 2007 Bennett \cite{Bennett07} sharpened this result employing several techniques such as the hypergeometric method of Thue and Siegel as well as a number of gap principles and proved that there are no $k$-th power Diophantine triple of the form $\{1,a,b\}$ for $k\geq 3$. Following these techniques, we prove the following theorem.
  \begin{theorem}\label{thm: main}
    Let $k\geq 3$ and $1<a<b<c$ be integers with $a^k<b<c$. Then there exists no $k$-th power Diophantine triple of the form $\{a^k, b, c\}$.
  \end{theorem}
\subsection*{Overview}\label{sec: Overview}
    The proof strategy follows closely the approach presented in \cite{Bennett07} and begins with the assumption that there exists a $k$-th power Diophantine triple as defined in Theorem~\ref{thm: main}. Next, we show that $c$ must be much larger than $a$ and $b$ and derive a good rational approximation of 
    $$
        \sqrt[k]{1+\frac{1}{a^k b}}.
    $$
    Then Lemma~\ref{lem: lower bounds} yields a contradiction for $k\geq 5$. As outlined in Bennett \cite{Bennett07}, for the cases $k=4$ and $k=3$, we introduce gap principles and show that $a^2t/s$ is a convergent of the infinite simple continued fraction expansion of $\sqrt[k]{a^kb}$. This yields a contradiction in all but a finite number of cases. 
    The remaining cases are then resolved using computer calculations of the continued fraction expansion of $\sqrt[k]{a^kb}$ with Sagemath \cite{Sage93}.

    \begin{remark}
        We note that for $a\geq 2$ setting $r:=a^k+1$ yields
        $$
            r^k = (a^k + 1 )^k = \sum\limits_{i=0}^{k}\binom{k}{i}a^{ik} = \sum\limits_{i=1}^{k}\binom{k}{i}a^{ik}+1 = a^k\sum\limits_{i=0}^{k-1}\binom{k}{i+1}a^{ik}+1.
        $$
        Also it holds that
        $$
            \sum\limits_{i=0}^{k-1}\binom{k}{i+1}a^{ik}> \binom{k}{2}a^k\geq a^k.
        $$
        Therefore, for all $a\geq 2$ setting $b:=\sum_{i=0}^{k-1}\binom{k}{i+1}a^{ik}$ yields a $k$-th power Diophantine pair $\{a^k,b\}$ with $1<a^k < b$, warranting our investigation into triples of the form $\{a^k,b,c\}$.
    \end{remark}

\section{Preliminaries}\label{sec: Prelim}

  In this section, we compile some basic statements and results that we will use later in the proof of our main theorem. 
  For $u > v >1$ and $k \geq 3$ it holds that
  \begin{equation}\label{eq: u^k-v^k}
    u^k-v^k = (u-v)\sum\limits_{i=0}^{k-1}u^{k-1-i}v^i > k(u-v).
  \end{equation}

      We will use the following result of Bennett~\cite[Theorem 1.3]{Bennett97}. For $n\geq 2$, define
      $$
        \mu_n = \prod\limits_{\substack{p\mid n\\ p \text{ prime}}} p^{1/(p-1)}.
      $$
      \begin{lemma}\label{lem: lower bounds}
        Let $p$ and $q$ be positive integers. If $n$ and $N$ are positive integers with $n\geq 3$ and 
        $$
          (\sqrt{N}+\sqrt{N+1})^{2(n-2)} > (n\mu_n)^n,
        $$
        then
        $$
          \left|\sqrt[n]{1 + \frac1N} - \frac{p}{q}\right| > (8n\mu_nN)^{-1}q^{-\lambda}
        $$
        with
        $$
          \lambda = 1 + \frac{\log(n\mu_n(\sqrt{N}+\sqrt{N+1})^2)}{\log((1/n\mu_n)(\sqrt{N}+\sqrt{N+1})^2)}.
        $$
      \end{lemma}

      Let $\alpha$ be an irrational number. We will write
      $$
        \alpha = [a_0, a_1,\dots ]
      $$
      for the infinite simple continued fraction expansion of $\alpha$.
      Here, $a_j$ denotes the $j$-th partial quotient satisfying $a_0\in \Z, a_j\in\N$ for $j\geq 1$.
      Furthermore, we write
      $$
        \frac{p_j}{q_j} = [a_0, a_1,\dots, a_j]
      $$
      for the $j$-th convergent of $\alpha$.
      We will use the following result on continued fractions; see e.g. \cite[Chapter 6]{Baker84}.
      \begin{lemma}
        Let $\alpha = [a_0, a_1,\dots].$ Then for $j\geq 1$ it holds that
        \begin{equation}\label{lem: continued fraction lower bound}
          \left|\alpha-\frac{p_j}{q_j}\right| > \frac{1}{(a_{j+1}+2)q_j^2}.
        \end{equation}
      \end{lemma}

\section{Proof of Theorem~\ref{thm: main}}\label{sec: Proof}

    Let us suppose that there exists a $k$-th power Diophantine triple of the form $\{a^k, b, c\}$ as described in Theorem~\ref{thm: main}, i.e. there exist non-zero integers $1<r<s<t$ such that
    \begin{equation}\label{eq: Diophantine triple}
      a^kb + 1 = r^k, \qquad  a^kc + 1 = s^k, \qquad bc + 1 = t^k.
    \end{equation}
    As for composite $k= pq$ the above equations may be written in the form
    $$
      (a^p)^qb + 1 = (r^p)^q, \qquad (a^p)^qc + 1 = (s^p)^q, \qquad bc + 1 = (t^p)^q,
    $$
    it suffices to consider the cases $k = 4$ and $ k\geq 3$ a prime in order to prove Theorem~\ref{thm: main}.
    We begin by establishing a lower bound for $a^kc$.
    
    Let now $k=4$ or $k\geq 3$ be a prime. From (\ref{eq: Diophantine triple}) it follows that
    $$
      r^k s^k = (a^k b+1)(a^k c+1)> a^{2 k}(b c+1)= a^{2 k} t^k.
    $$
    From this we obtain $rs > a^2 t$, which implies that there exists a non-zero positive integer $z$ such that
    \begin{equation}\label{eq: rk = a^2 t+z }
      rs = a^2t+z.
    \end{equation}
    This leads to the equality
    \begin{align}\label{eq: binomial expansion}
        r^k s^k & = (a^2t+z)^k = \left( (a^{2k}(bc + 1))^{\frac1k} + z \right)^k\\
        & = a^{2k}(bc + 1) + k(a^{2k}(bc + 1))^{\frac{k-1}{k}}z \nonumber \\
        &\qquad +\binom{k}{2}(a^{2k}(bc + 1))^{\frac{k-2}{k}}z^2 + \cdots + z^k.\nonumber
    \end{align}
    We have
    \begin{align}
      r^k s^k &= (a^k b + 1)(a^k c + 1) = a^{2 k} b c + a^k b + a^k c + 1 \nonumber \\
      & \leq a^{2 k} (b c + 1) + a^k b + a^k c + 1 \nonumber
    \end{align}
    and, thus, it follows from (\ref{eq: binomial expansion}) that 
    \begin{align}\label{eq: binomial inequality}
        a^k b+a^k c+1 &\geq k(a^{2 k}(b c+1))^{\frac{k-1}{k}}z+\binom{k}{2}(a^{2 k}(b c+1))^{\frac{k-2}{k}}z^2 + \cdots \\
        &\qquad \cdots + z^k. \nonumber
    \end{align}
    Assume that $a^k b \geq(a^{2 k}(b c + 1))^{\frac{k-2}{k}}$. Then for $k\geq 4$ we obtain
    $$
      (a^k b)^k \geq (a^{2 k}(b c + 1))^{k-2} \geq a^{2k(k-2)} (bc)^{k-2},
    $$
    which implies
    $$
      b^2 \geq a^{k(k-4)}c^{k-2} > c ^2
    $$
    contradicting $ b < c$.
    For the case $ k = 3 $, we get the inequality
    $$
       (a^3b)^3\geq (a^6(bc+1))>a^6bc,
    $$
   which yields
    $$
      a^3 b^2 > c.
    $$
    From (\ref{eq: binomial inequality}) it follows that
    $$
      3(a^{6}bc)^{\frac23}z <3(a^{6}(bc+1))^{\frac23}z\leq a^3b+a^3c+1\leq 3a^3c,
    $$
    which implies $ a^3 b^2 \leq c$, a contradiction.
    Thus, for $k\geq 3$ we have
    \begin{equation}\label{eq: auxiliary a^k b inequality}
      a^k b < (a^{2k}(b c+1))^{\frac{k-2}{k}}.
    \end{equation}
    Now, assume that $a^k c \leq k (a^{2 k}(b c+1))^{\frac{k-1}{k}} z$.
    Then, from the inequality (\ref{eq: binomial inequality}) we obtain
    \begin{align}
      a^k b+a^k c+1 &\geq k(a^{2 k}(b c+1))^{\frac{k-1}{k}}z+\binom{k}{2}(a^{2 k}(b c+1))^{\frac{k-2}{k}}z^2 + \cdots \nonumber \\
        &\qquad\cdots + k(a^{2k}(bc+1))^{\frac{1}{k}}z^{k-1} + z^k\nonumber \\
        & \geq a^k c + (a^{2k}(bc+1))^{\frac{k-2}{k}} +1,\nonumber
    \end{align}
    which leads with (\ref{eq: auxiliary a^k b inequality}) to the contradiction
    $$
      a^k b \geq (a^{2k}(bc+1))^{\frac{k-2}{k}} > a^kb.
    $$
    Therefore, it holds that $a^k c > k (a^{2 k}(b c+1))^{\frac{k-1}{k}} z$.
    This implies
    $$
      a^k c > k(a^{2k}(bc+1))^{\frac{k-1}{k}}z > k(a^{2k} b c)^{\frac{k-1}{k}}z.
    $$
    Finally, we can conclude
    \begin{equation}\label{eq: lower bound a^k c}
      a^k c > k^k z^k (a^k b)^{k-1}.
    \end{equation}
    Next, we derive an upper bound for $a^k c$, which alongside (\ref{eq: lower bound a^k c}) will lead to a contradiction.
    From (\ref{eq: Diophantine triple}) it follows that
    \begin{align}
      0 < \frac{r^k s^k}{(a^2t)^k} - \frac{a^k b+1}{a^k b} &=\frac{(a^k b+1)(a^k c+1)}{a^{2k}(b c+1)}-\frac{a^k b+1}{a^k b} \nonumber \\
        & =\frac{a^k b^2+b-a^{2k} b-a^k}{a^{2k}b(b c+1)} \nonumber \\
        & < \frac{a^k b}{a^{2k}(bc+1)} \nonumber \\
        & =\frac{a^kb}{(a^2 t)^k}.\nonumber
    \end{align}
    This, combined with (\ref{eq: u^k-v^k}), yields
    \begin{equation}\label{eq: upper bound dio. approximation}
       \left|\sqrt[k]{1+\frac{1}{a^k b}} -\frac{r s}{a^{2}t}\right|<\frac{a^k b}{k (a^2t)^k}.
    \end{equation}
    We want to employ Lemma~\ref{lem: lower bounds} with $N= a^kb$ and $q=a^2t$, but first we need to check for which cases the condition
    \begin{equation}\label{eq: dio. approximation condition}
      (\sqrt{N}+\sqrt{N+1})^{2(k-2)} > (k\mu_k)^k
    \end{equation}
    are satisfied. We have $N= r^k -1$ and $r^k> a^k b \geq 2^k(2^k+1)>2^{2k}$, which means that $r\geq 5$. Thus, it suffices to ensure that
    $$
      (N+1)^{(k-2)} = r^{k(k-2)}\geq 5^{k(k-2)}>(k\mu_k)^k
    $$
    holds to guarantee (\ref{eq: dio. approximation condition}).
    This inequality is satisfied for $k = 4$ and $k\geq 5$ prime.
    For $k=3$ and $r = 5$ one can directly compute that (\ref{eq: dio. approximation condition}) holds.
    Therefore, we may apply Lemma~\ref{lem: lower bounds} in all of our cases and obtain in combination with (\ref{eq: upper bound dio. approximation})
    $$
      (8k\mu_k a^kb)^{-1}(a^2t)^{-\lambda}< \frac{a^k b}{k(a^2 t)^k}.
    $$
    This simplifies to
    \begin{equation}\label{eq: a^2 t upper bound}
      (a^2 t)^{k -\lambda} < 8 \mu_k (a^k b)^2
    \end{equation}
    and finally
    \begin{equation}\label{eq: a^kc upper bound}
      (a^k c)^{k-\lambda} < 8^k\mu_k^k (a^k b)^{k+\lambda}.
    \end{equation}
    In the following subsections, we will differentiate between the cases $k =3,4$ and $k\geq5$ a prime.
    
  \subsection{Case $ k\geq 5 $}\label{sec: k >= 5}
    For a prime $k$ it holds that $\mu_k = k^{\frac{1}{k-1}}$.
    For $N=r^k-1=a^kb$, $\lambda$ is monotonously decreasing in $k$ and $r$. Now, let $k\geq 7$ be a prime.
    For $r\geq 5$, we have $\lambda < 2.44$.
    Thus, (\ref{eq: a^kc upper bound}) in combination with (\ref{eq: lower bound a^k c}) yields
    $$
      (a^k b)^{k-1} < 8^{\frac{k }{k-2.44}}(a^k b)^{\frac{k+2.44 }{k-2.44}}< 28 (a^k b)^{2.1}.
    $$
    This leads to a contradiction for $k\geq 7$ and $r\geq 5$.
    If $k=5$ and $r\geq 5$, then $\lambda < 2.55$, which, as above, leads to
    $$
      (a^5 b)^{4}  < 79 (a^5 b)^{3.1}.
    $$
    This is a contradiction for $r\geq 5$ and therefore we have proven Theorem~\ref{thm: main} for $k \geq 5$ a prime.
  
  \subsection{Case $ k = 4 $}\label{sec: k = 4}
    We now consider the case $k=4$.
    As inequality (\ref{eq: lower bound a^k c}) is not strong enough on its own to lead to a contradiction, we introduce a gap principle to sharpen this inequality.
    First, assume that $s\leq 3 r^3 z$. Then from (\ref{eq: lower bound a^k c}) it follows that
    $$
      4^4 z^4 (a^4 b)^3 < a^4c < s^4 \leq 3^4(r^4)^3 z^4 = 3^4 (a^4b+1)^3z^4,
    $$
    which implies
    $$
       \left(\frac{4}{3}\right)^4 < \left(\frac{a^4b+1}{a^4b}\right)^3.
    $$
    This is satisfied for $a^4b\leq 2$, which contradicts $1 < a < b$. Therefore, we have
    \begin{equation}\label{eq: k = 4 s lower bound}
      s > 3 r^3 z.
    \end{equation}
    From (\ref{eq: rk = a^2 t+z }) it follows that
    $$
      (r^4-1)(s^4-1)=a^8 (t^4-1) =  (a^2 t)^4 - a^8 = (r s-z)^4-a^8,
    $$
    and thus
    \begin{equation}\label{eq: k = 4 rsz equality}
      4(r s)^3 z + 4 r s z^3 + 1 + a^8 = r^4 + s^4 + 6(r s z)^2 + z^4.
    \end{equation}
    Consequently, we have
    \begin{equation}\label{eq: k = 4 rsz mod s^3}
      4 r s z^3 + 1 + a^8 \equiv r^4 + 6(r s z)^2 + z^4\mod{s^3}.
    \end{equation}
    From (\ref{eq: k = 4 rsz mod s^3}) we have the following three cases, and in each case we derive $z > r/2$. 
    First, in the case where
    $$
      4 r s z^3 + 1 + a^8 = r^4 + 6(r s z)^2 + z^4,
    $$
    it follows from (\ref{eq: k = 4 rsz equality}) that $ s^4 = 4(rs)^3z$ and therefore $ s = 4 r^3 z$.
    This implies that $ r^4 \mid z^4 - 1 - a^8 $ and hence $ z > r/2$.
    Second, if 
    $$
      4 r s z^3 + 1 + a^8 \geq r^4 + 6(r s z)^2 + z^4 + s^3,
    $$
    we can conclude
    \begin{align}
      4 r s z^3 + 1 + a^8 &> 1 + 6 (a^4 b+1)^{\frac{1}{2}}(a^4 c+1)^{\frac{1}{2}} + s^3 \nonumber \\
      &> 1 + (a^4 b+1) + s^3 \nonumber \\
      &> 1 + a^8 + s^3. \nonumber
    \end{align}
    This together with (\ref{eq: k = 4 s lower bound}) leads to
    $$
      4 r z^3 > s^2 > 3^2 r^6 z^2,
    $$
    which yields $ z > r/2 $.
    Third, in the case where
    $$
      4 r s z^3 + 1 + a^8 + s^3 \leq r^4 + 6(r s z)^2 + z^4 + s^3,
    $$
    suppose that $s^3 \geq 6(rsz)^2$.
    This implies
    $$
       4 r s z^3 + 1 + a^8 + s^3 \leq r^4 + 6(r s z)^2 + z^4 \leq r^4 + s^3 + z^4
    $$
    and in combination with (\ref{eq: k = 4 s lower bound}) yields the contradiction
    $$
      r^4 + z^4 < 12 r^4 z^4< 12 r^4 z^4 + 1 + a^8 < 4 r s z^3 + 1 + a^8 \leq r^4 +z^4.
    $$
    Therefore, it holds that $ s^3 < 6(rsz)^2 $ and (\ref{eq: k = 4 s lower bound}) allows us to conclude that $z > r/2$ in all cases.
    Consequently, by (\ref{eq: lower bound a^k c}) we have
    $$
       a^4 c  > 4^4 z^4(a^4 b)^3 > 16 r^4 (a^4 b)^3 > 16 (a^4 b)^4.
    $$
    We note that $\mu_4 = 2$ and combine the above inequality with (\ref{eq: a^kc upper bound}) to obtain
    $$
      (r^4-1)^{12-5\lambda} < 16^\lambda,
    $$
    which for $r\geq 35$, and thus $\lambda < 2.308$, leads to a contradiction.

    Now we consider $5\leq k \leq 34$. From Lemma~\ref{lem: lower bounds} it follows that
    $$
      (a^2 t)^{4-\lambda} < 16 (a^4b)^2
    $$
    and thus
    $$
      a^2 t < \left(16(r^4-1)\right)^{\frac{1}{4-\lambda}}.
    $$
    Furthermore, for $5 \leq r \leq 34$, we have
    \begin{equation}\label{eq: k = 4 10^8 bound}
      \left(16(r^4-1)\right)^{\frac{1}{4-\lambda}} < 10^8.
    \end{equation}
    Additionally, (\ref{eq: Diophantine triple}) implies for $k\geq 3$
    $$
      (r^k-1) s^k - a^{2k} t^k = a^k b - a^{2k} < a^k b - 1 = r^k-2
    $$
    and we have
    \begin{equation}\label{eq: k = 4 dio approx basic}
      (r^k-1)-\frac{a^{2k} t^k}{s^k} < \frac{r^k-2}{s^k}. 
    \end{equation}
    By (\ref{eq: u^k-v^k}) it holds that
    \begin{align}\label{eq: k = 4 u^k-v^k expansion}
        (r^k-1)-\frac{a^{2k} t^k}{s^k} = &\left(\sqrt[k]{r^k-1}-\frac{a^2 t}{s}\right)\left( (r^k-1)^{\frac{k-1}{k}}\right.  \\
        &\qquad \left. + (r^k-1)^{\frac{k-2}{k}}\frac{a^2 t}{s} + \cdots + \frac{a^{2(k-1)} t^{k-1}}{s^{k-1}}\right). \nonumber
    \end{align}
    Moreover, we have the inequality
    $$
        s^k (r^k-2)=(a^k c + 1)(a^k b - 1) < a^{2k} b c + a^{2k} = a^{2k} (b c+1) = a^{2k} t^k         
    $$
    and thus
    \begin{equation}\label{eq: k = 4 extra inequality}
        \frac{a^2 t}{s}>\sqrt[k]{r^k-2}.
    \end{equation}
    Therefore, by combining (\ref{eq: k = 4 dio approx basic}), (\ref{eq: k = 4 u^k-v^k expansion}) and (\ref{eq: k = 4 extra inequality}) we obtain
    $$
        \frac{r^k-2}{s^k} > (r^k-1)-\frac{a^{2k} t^k}{s^k} > \left(\sqrt[k]{r^k-1}-\frac{a^2 t}{s}\right)k(r^k-2)^{\frac{k-1}{k}},
    $$
    which leads to 
    \begin{equation}\label{eq: k = 4 dioph approx upper bound}
        \left| \sqrt[k]{r^k-1}-\frac{a^2 t}{s}\right|< \frac{(r^k-2)^{\frac{1}{k}}}{k\,s^k}.
    \end{equation}
    Then, in the case $k = 4$, from (\ref{eq: k = 4 dioph approx upper bound}) and (\ref{eq: k = 4 s lower bound}) we have
    $$
       \left| \sqrt[4]{r^4-1}-\frac{a^2 t}{s}\right| < \frac{1}{2s ^2}.
    $$
    Thus, $a^2 t / s = p_j/q_j$ is a convergent in the continued fraction expansion of $\sqrt[4]{r^4-1}$.
    Furthermore, it holds that
    $$
      \sqrt[4]{r^4-1} > \frac{a^2 t}{s},
    $$
    which means that $j$ is even.
    By Lemma~\ref{lem: continued fraction lower bound} we have
    $$
      \left| \sqrt[4]{r^4-1}-\frac{a^2 t}{s}\right| > \frac{1}{(a_{j+1}+2)s^2},
    $$
    which by (\ref{eq: k = 4 s lower bound}) in combination with (\ref{eq: k = 4 dioph approx upper bound}) leads to
    $$
      9 r^8 < 4 s^4 < (a_{j+1} + 2) (r^4-2)^{\frac14} < (a_{j+1} + 2) r.
    $$
    Therefore, the inequality
    \begin{equation}\label{eq: k = 4 quotient condition}
      a_{j+1} > 9 r^7 - 2 \geq 703123      
    \end{equation}
    holds.
    For all $5\leq r \leq 35$ we can calculate that $ p_{13} > 10^8$ and, thus, according to (\ref{eq: k = 4 10^8 bound}) we have $j\in\{0, 2, 4, 6, 8, 10, 12\}$.
    However, for no such $j$ is inequality (\ref{eq: k = 4 quotient condition}) satisfied, which is a contradiction.
    Therefore, we have proven Theorem~\ref{thm: main} for the case $ k = 4$.

  \subsection{Case $ k = 3 $}\label{sec: k = 3}

    Finally, we consider the case $k = 3$.
    We begin similarly as in the previous case to sharpen the inequality (\ref{eq: lower bound a^k c}).
    By (\ref{eq: k = 4 dioph approx upper bound}) we have
    \begin{equation}\label{eq: k = 3 dioph approx upper bound}
      \left|\sqrt[3]{r^3 - 1}-\frac{a^2 t}{s}\right|<\frac{\left(r^3-2\right)^{\frac{1}{3}}}{3 s^3}.
    \end{equation}
    From (\ref{eq: lower bound a^k c}) we have $ s > 3 (a^3 b)^{\frac23} = 3 (r^3-1)^{\frac23}$ and thus
    $$
      \left|\sqrt[3]{r^3-1}-\frac{a^2 t}{s}\right| < \frac{1}{2s^2}
    $$
    Therefore, $a^2 t /s = p_j/q_j$ must be a convergent of the continued fraction expansion of $\sqrt[3]{r^3-1}$.
    Additionally, as
    $$
      \sqrt[3]{r^3-1} > \frac{a^2 t}{s}
    $$
    the index $j$ is even.
    By Lemma~\ref{lem: continued fraction lower bound} it holds that
    $$
      \left|\sqrt[3]{r^3-1}-\frac{a^2 t}{s}\right| > \frac{1}{(a_{j+1} + 2)s^2}.
    $$
    Consequently, combining the above inequality with (\ref{eq: k = 3 dioph approx upper bound}) and (\ref{eq: lower bound a^k c}) we obtain
    $$
      9(r^3 - 1 )^{\frac23} < 3 s < (a_{j+1} + 2 ) (r^3 - 2)^{\frac13}
    $$
    and thus
    \begin{equation}\label{eq: k = 3 quotient condition}
      a_{j+1} > 3 r -2 \geq 13.
    \end{equation}
    Now, we examine the continued fraction expansion of $\sqrt[3]{r^3-1}$ more closely.
    In Bennett~\cite[Chapter 4]{Bennett07} the first partial quotients for $r\geq 3$ are computed as
    $$
      a_0 = r - 1, \quad a_1 = 1, \quad a_2 = 3 r^2 - 2, \quad a_3 = 1, \quad a_4 = r - 2, \quad a_5 = 1.
    $$
    As $j$ is even and (\ref{eq: k = 3 quotient condition}) holds, we have $j\geq 6$.
    Further, analysis of the partial quotients, see~\cite[Chapter 4]{Bennett07}, in combination with (\ref{eq: k = 3 quotient condition}) leads to $j\geq8$ and thus $s\geq q_8 >5r^6$, except, possibly, for
    \begin{equation}\label{eq: k = 3 exception set}
      r \in\{2,3,5,7,9,11,15,17,19,21,25,27,31,37,41,47,57\}.
    \end{equation}
    In our case, it holds that $r\geq 5$.
    Furthermore, we have $r^3-1 = a^3b$ and quick computations confirm that there is no $r < 9$ of this form.
    Next, for the remaining values of $r$ in (\ref{eq: k = 3 exception set}) one can directly compute $a_9\leq 10$.
    Therefore, (\ref{eq: k = 3 quotient condition}) again guarantees that
    $$
      s \geq q_{10} > q_8 > 5r^6
    $$
    for $r\geq 9$.
    From this we can derive the inequality
    $$
      a^3 c > 125 (a^3 b)^6,
    $$
    which, by applying (\ref{eq: lower bound a^k c}), results in
    $$
      125^{3-\lambda} (a^3 b)^{18 - 6\lambda} < 8^3 3^{\frac32} (a^3 b)^{3+\lambda}.
    $$
    Thus, we have
    $$
      (a^3 b)^{15 - 7\lambda} < 8^3 3^{\frac32} 125^{\lambda - 3},
    $$
    which is contradicted by $r\geq 7973$.
    For $9\leq r\leq 7972$, from (\ref{eq: a^2 t upper bound}) it follows that
    $$
       a^2 t < \left(8 \sqrt{3} (a^3 b)^2\right)^{\frac{1}{3-\lambda}}
    $$
    and for $a^3 b = r^3-1$ and $\lambda$ as defined in Lemma~\ref{lem: lower bounds} it holds that
    \begin{equation}\label{eq: k = 3 p_j upper bound}
      p_j = a^2 t < 10^{32}.
    \end{equation}
    Additionally, it holds that
    $$
      (r^3-1)q_j^3 - p_j^3 = a^3 b s^3 - a^6 t^3  =   r^3 - 1 -a^6,
    $$
    which implies
    \begin{equation}\label{eq: k = 3 q_j p_j divisibility condition}
      (r^3-1)q_j^3 - p_j^3 \mid r^3 - 1 -a^6.
    \end{equation}
    Now, we can systematically eliminate all remaining candidates for $r$.
    Simple calculations show that for $9\leq r\leq 7972$ only $1892$ are of the form $ r^3 -1 = a^3 b$.
    For the remaining cases, we can verify with inequality (\ref{eq: k = 3 p_j upper bound}) that none satisfies the condition (\ref{eq: k = 3 q_j p_j divisibility condition}).
    Therefore, we have proven Theorem~\ref{thm: main} for $k\geq 3$. \hfill $\square$

\Acknowledgements{The authors would like to thank the Bundesministerium f\"{u}r Frauen, Wissenschaft und Forschung (BMFWF) for the support under the Sparkling Science project SPA 01-080 `MAJA -- Mathematische Algorithmen für Jedermann Analysiert'.}

\end{document}